\documentclass[12pt]{amsart}
\usepackage{graphicx}
\usepackage{tabularx}
\usepackage{ragged2e}
\usepackage{url}
\usepackage{hyperref}
\usepackage{amsmath,color}
\usepackage{fullpage}
\usepackage{microtype}
\newtheorem{theorem}{Theorem}[section]
\newcommand{\seqnum}[1]{\href{http://oeis.org/#1}{\underline{#1}}}

\usepackage[greek,english]{babel}
\newcommand{\numbers}{\ensuremath{\alpha\rho\iota\theta\mu o\varsigma}}
\newcommand{\magnitude}{\ensuremath{\mu\varepsilon\gamma\varepsilon\theta o\varsigma}}
\begin{document}

\title{What is the smallest prime?}%

\author{Chris K.~Caldwell}%
\address{University of Tennessee at Martin}%
\email{caldwell@utm.edu}%

\author{Yeng Xiong}%
\address{University of Tennessee at Martin}%
\email{yenxion@ut.utm.edu}%

\begin{abstract}
   What is the first prime?  It seems that the number two should be the obvious answer, and today it is, but it was not always so.  There were times when and mathematicians for whom the numbers one and three were acceptable answers.  To find the first prime, we must also know what the first positive integer is.  Surprisingly, with the definitions used at various times throughout history, one was often not the first positive integer (some started with two, and a few with three).  In this article, we survey the history of the primality of one, from the ancient Greeks to modern times. We will discuss some of the reasons definitions changed, and provide several examples.  We will also discuss the last significant mathematicians to list the number one as prime.
\end{abstract}

\maketitle

\section{Introduction}\label{sect intro}

In our research, we kept running across claims that the number one ``used to be prime.''  For example, at the time we wrote this article, the sequence \seqnum{A008578} in OEIS \cite{OEIS2012} was named ``prime numbers at the beginning of the 20th century (today 1 is no longer regarded as a prime).''  And before we changed it in late 2011, the English Wikipedia article on \textit{prime number} \cite{OldWiki} said % previous version Revision as of 02:26, 4 September 2011 (edit) Caldwell's edit: Revision as of 15:15, 8 September 2011 from 208.87.72.31
``until the 19th century, most mathematicians considered the number 1 a prime.''  There were (and are) individuals who defined the number one to be prime, but as we will see in this brief note, the history is much more complicated than ``one used to be prime.''  At first, most started the sequence of primes with 2, but some began with 3 and a rare few began with 1.  There does not appear to be any period of time during which most mathematicians defined one to be a prime.

Before we start, let's make two points clear. First, whether or not a number (especially unity) is a prime is a matter of definition, so a matter of choice, context and tradition,  not a matter of proof.  Yet definitions are not made at random; these choices are bound by our usage of mathematics and, especially in this case, by our notation.  We will see, for example, that as the use of the integers changed, one (unity, the monad), which at first was not generally viewed as a number by the ancient Greeks, `becomes' a number some 2,000 years later.  For much of history it did not even make sense to ask if the number one was a prime.

Second, when we stick to the (modern) positive integers, it is not difficult to adjust our theorems to allow the number one to be called a prime.  For example, some say unity is omitted from the primes to preserve the uniqueness of factorization.  However, in the past, mathematicians who defined one to be prime just added a couple words to the fundamental theorem of arithmetic (FTA).  If we start with the version of the FTA stated by Crandall and Pomerance \cite[p.~1--2]{CP2005}), we need only add the two characters ``$1 <$'':
\begin{theorem} For each natural number $n$ there is a unique factorization $$n = p_1^{a_1}p_2^{a_2} \cdots p_k^{a_k},$$ where exponents $a_i$ are positive integers and $\textcolor{red}{1 <}\,\, p_1 < p_2 < \cdots < p_k$ are primes. \end{theorem}
\noindent The real problem with one being a prime only became apparent in the 18th century when mathematicians expanded their purview to wider notions of integers, such as the Gaussian integers and general number fields.  In that context, they were forced to address the notion of units, to separate the concepts of irreducibility and primality, etc. The generalization of prime to unique factorization domains clarified the role of unity and now informs the way we define primality in the ordinary positive integers.

\section{Before 1900}\label{sect pre}

To get started, we will show that in the early days of defining primes, which we will illustrate with Euclid, the unit (one) was not even considered a number, so was \textit{ipso facto} excluded from the primes.  This begins to change when Stevin showed how to calculate with decimal real numbers in 1585, bringing on a confused period in which some held unity to be a prime and others did not. When Gauss first states and proves the fundamental theorem of arithmetic (FTA) in 1801, writers begin to coalesce slowly around the modern definition of prime---which again excludes unity.

\subsection{Before the number one}

Euclid, for example, defines unity ($\mu o\nu\acute{\alpha}\delta\iota$), numbers (\numbers), and the primes ($\pi\rho\tilde{\omega}\tau o\varsigma$) in book 7 of his \textit{Elements} as follows:
\begin{itemize}
   \item A unit is that by virtue of which each of the things that exist is called one.
   \item A number is a multitude composed of units.
   \item A prime number is that which is measured by a unit alone.
\end{itemize}
Euclid (and his contemporaries) did not need to say explicitly that one was not a prime because primes were a subcategory of the numbers, so one was not a number.  Smith \cite[p.~20]{Smith1958} states it this way:
\begin{quote}
   Aristotle, Euclid, and Theon of Smyrna defined a prime number as a number ``measured by no number but an unit alone,'' with slight variations of wording.  Since unity was not considered a number, it was frequently not mentioned.
\end{quote}
Writers at times broke this silence; for example Martianus Capella (c.400, \cite[pp.~285--286]{Stahl1992}), an early and influential developer of the system of seven liberal arts, summarizes the properties of the first integers as follows:
\begin{quote}
   We have briefly discussed the numbers comprising the first series, the deities assigned to them, and the virtue of each number. I shall now briefly indicate the nature of number itself, what relations numbers bear to each other, and what forms they represent. [\ldots] Let us consider all numbers of the first series according to the above classifications: the monad is not a number; the dyad is an even number; the triad is a prime number, both in order and in properties; the tetrad belongs in the even times even class; the pentad is prime; the hexad belongs to the odd times even or even times odd (so it is called perfect); the heptad is prime; the octad belongs to [\ldots]
\end{quote}
Martinus and others such as Nicomachus (c.100) and Iamblichus (c.300) \cite[p.~73]{Heath1981}, Boethius (c.500) \cite[pp.~89--95]{Masi1983}, and Cassiodorus (c.550) \cite[p.~5]{Grant1974} make the primes a subset of the odds, excluding both one and two from the primes---so for them, the smallest prime was three.  (It is easy to extend this list of those for whom the first prime was three well into the 16th century \cite{CRX2012}.) Most of the ancient Greeks however, like Euclid, began the sequence of primes with two.

\begin{figure}[hbt]
  \centering
  \includegraphics{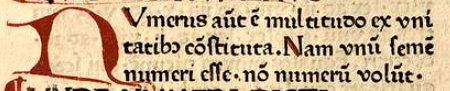}
  \caption{A 15th century copy of Isidore of Seville's \textit{Etymologiae} (636)}
  \label{fig_Isadorus}
\end{figure}

For about 2000 years most writers varied little from the position of the ancient Greeks.  For example, Figure~\ref{fig_Isadorus} shows a copy of Isidore of Seville's (c.636) work made in 1493.
Here Isidore (\cite[pp.~4--5]{Grant1974}) writes ``Number is a multitude made up of units. For one is the seed of number but not number.''   Nine centuries after Isidore,
J.~K\"obel (1537) repeats this argument ``Darauss verstehstu das 1. kein zal ist / sonder es ist ein gebererin / anfang / und fundament aller anderer zalen''
(see Menninger \cite[p.~20]{Menninger1992}). Simply put, for them numbers are multiples and therefore divisible, but unity is neither a multiple or divisible.  (If you would like a longer argument, al-Kind{\=\i} (c.850) argues at great length \cite{AlKindi1974} that one is not a number, so that it is neither even nor odd\ldots.)

Omitting one from the numbers, hence the primes, was the general rule, but there were a few isolated exceptions.   Speusippus (c.350BC), Tar\'an writes \cite[p.~276]{Taran1981},
\begin{quote}
   is exceptional among pre-Hellenistic thinkers in that he considers one to be the first prime number. And Heath, \textit{Hist. Gr. Math.}, I, pp.~\oldstylenums{69}-\oldstylenums{70}, followed by Ross, \textit{Aristotle's Physics}, p.~\oldstylenums{604}, and others, is mistaken when he contends that Chrysippus, who is said to have defined one as $\pi\lambda\tilde{\eta}\theta o\varsigma$ \textgreek{>'en} (cf. Iamblichus, \textit{In Nicom. Introd. Arith.}, p.~II, \oldstylenums{8}-\oldstylenums{9} [Pistelli]), was the first to treat one as a number (cf. further p.~\oldstylenums{38}f. with note \oldstylenums{189} \textit{supra}).
\end{quote}
A later exception may be Rabbi ben Ezra (c.1140) in his \textit{Sefer ha-Echad} \cite[pp.\ 27--28]{Smith1958}.

\subsection{One becomes a number}

\begin{figure}[hbt]
  \centering
  \includegraphics{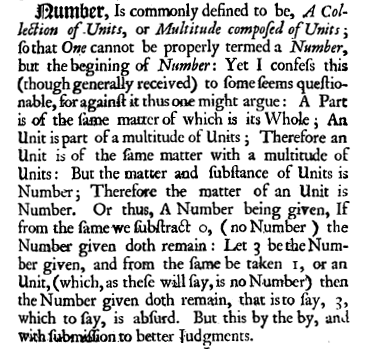}
  \caption{Moxon's definiton of number (1679)}
  \label{fig_moxon}
\end{figure}

In 1585 Simon Stevin published his tract \textit{De Thiende} which laid the basis for modern decimal expansions \cite{KK2011}.  Here Stevin showed how to represent and operate on both the positive integers (Euclid's \numbers) and the magnitudes (lengths of geometrical objects, Euclid's \magnitude) with the same notation and algorithms.  In his thesis ``The Concept of \underline{One} as a Number,'' C.~Jones \cite[p.~299]{Jones1978} put it this way:
\begin{quote}
   In general, mathematics before Stevin is of one character and, after him, it is of another reflecting his contributions.  In this regard, he is like Euclid: he stood at a watershed in the history of mathematics.  And as with Euclid, he was so successful that, from our present day vantage point, it is hard to see the other side of that watershed.  Over there, \underline{one} was not a number; here and now, it is; even $\pi$ is a number, and $i$, and aleph null.
\end{quote}
No change happens instantly, and nearly a century later when Moxon wrote the first English mathematical dictionary (\cite[p.~97]{Moxon1679}), Moxon defined number (with respect to unity) with trepidation (see Figure~\ref{fig_moxon}).
%\begin{quote}
%   {\vold {\gothfamily \large Number}, Is= commonly defined to be, \textit{A Collection of Units=}, or \textit{Multitude composed of Units=}; so that \textit{One} cannot be properly termed a \textit{Number}, but the beginning of \textit{Number}: Yet I confess= this= (though generally received) to some seems= questionable, for against it thus= one might argue: A Part is= of the same matter of which is= its= Whole; An Unit is= part of a multitude of Units=; Therefore an Unit is= of the same matter with a multitude of Units=: But the matter and substance of Units= is= Number; Therefore the matter of an Unit is= Number.  Or thus=, A Number being given, If from the same we subtract 0, (no Number) the Number given doth remain: Let 3 be the Number given, and from the same be taken 1, or an Unit, (which, as= these will say, is= no Number) then the Number given doth remain, that is= to say, 3, which to say, is absurd.  But this= by the by, and with submission to better Judgments=.} \\
%\end{quote}
The argument Moxon gives here, that 1 is a number otherwise $3-1=3,$ was one of Stevin's arguments, but it is best viewed as a metaphor for Stevin's most convincing message: notationally and algorithmically, there is no difference between one and the other numbers.

\begin{figure}[bth]
  \centering
  \includegraphics[width=.9\textwidth]{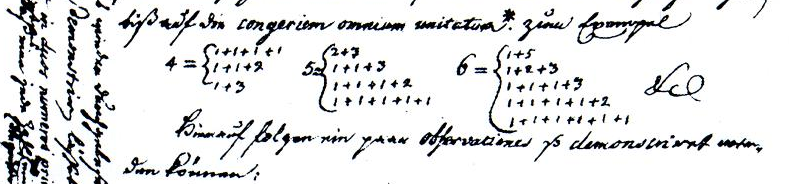}
  \caption{A part of Goldbach's letter to Euler writing integers as sums of primes}
  \label{fig_Goldbach}
\end{figure}

The redefinition of unity as a number now made it reasonable to ask if one was a prime; however, the way primes were used in this period did not force any particular definition.  In their history of the fundamental theorem of arithmetic, A{\u{g}}arg{\"u}n and {\"O}zkan write that at this time ``prime factorization was not looked upon as something of interest in its own right, but as a means of finding divisors'' \cite[p.~211]{AO2001}.  Similarly, A{\u{g}}arg{\"u}n and Fletcher, writing about the existence and uniqueness components of the FTA, note \cite{AF1997} the following:
\begin{quote}
   \ldots it is significant that Propositions VII.31 and VII.30 of the Elements lead immediately to their proofs although Euclid forbears to take these steps.  The first known proof of existence is due to al-Farisi (died c.1320), but he did not go on to prove uniqueness, mainly because his interest was in the divisors of a number rather than the factorisation itself \ldots And if a mathematician's interest is in greatest common divisors, or perfect numbers, or amicable numbers then the divisors are the crucial objects, whereas the prime factorisation is just a means to an end.
\end{quote}
Without clear direction on the primality of one, a period of confusion begins.  Many stuck with the tradition of excluding one from the primes, for example:
P.~A.~Cataldi (1603) \cite[p.~40]{Cataldi1603},
M.~Mersenne (1625) \cite[pp.~298--299]{Mersenne1625},
L\'eon de Saint-Jean (1657) \cite[p.~581]{Leon1657},
F.~v.~Schooten (1657) \cite[pp.~393--403]{Schooten1657},
the Shuli Jingyun (c.1720) \cite{Roegel2011} and
L.~Euler (1770) \cite[pp.~14--16]{Euler1770a}.
However, it was also reasonable to define the primes as the complement of the set of composites, often called the incomposites, such as in the influential tables of Brancker \& Pell (1688) \cite[p.~367]{Maseres1795}.  Others who listed one as prime in this period are
F.~Wallis (1685) \cite[p.~292]{Maseres1795},
J.~Prestet (1689) \cite[p.~141]{Prestet1689},
C.~Goldbach (1742) \cite{Goldbach1742} (see Figure~\ref{fig_Goldbach}),
J.~H.~Lambert (1770) \cite[p.~73]{Lambert1770},
A.~Felkel (1776) \cite{Felkel1776} and
E.~Waring (1782) \cite[p.~379]{Waring1782}.
We could make these lists much longer (see \cite{CRX2012} and \cite{RX2012}), but these are sufficient to illustrate the point.  Keep in mind that these authors may, or may not, have stuck to the indicated view of the primality of one; this list is just what they wrote at the indicated time.

Definitions depend on context (as well as on tradition and usage), so for example, the later number theorist V.\ A.\ Lebesgue (1791--1875) omitted unity from the primes in an 1856 article \cite{Lebesgue1856}, included it in an 1859 article \cite[p.~5]{Lebesgue1859}, and then excluded it in 1864 article \cite[p.~12]{Lebesgue1864}.  For many authors we must infer their position by their usage. For example, the snippet of Goldbach's famous letter to Euler in Figure~\ref{fig_Goldbach} about expressing integers as sums of prime makes his stand obvious, as does Gauss counting 168 primes below 1,000.

When we look back on this period, it seems odd that factorizations were used without asking if those factorizations were unique.  Our modern notions of integers, primes and unique factorization appear so obvious to us that we are blinded by hindsight.

\subsection{Factorization becomes central}

Gauss' \textit{Disquisitiones Arithmeticae} provides another watershed moment in the history of the integers.  A{\u{g}}arg{\"u}n and Fletcher write \cite{AF1997}
\begin{quote}
   It is not too much of an exaggeration to say that the result [FTA] passed from being unknown to being obvious without a proof passing through the head of any mathematician.  [\ldots] The first clear statement and proof of the FTA seems to have been given by Gauss (1775--1855) in his \textit{Disquisitiones Arithmeticae} of 1801.
\end{quote}
Gauss not only stated FTA (and gave it a partial proof), but while doing so laid the ground work for making unique factorization central to our understanding of integers.  Equally important, he discussed generalizations of the integers, such as the Gaussian integers, which made mathematicians look more carefully at the role of one (and all the units) in factorizations.  Gauss' text never defines primes explicitly, but still provides an eloquent argument for the modern definition of prime.

Even after Gauss' pivotal text, many continued to write that unity was prime.  Among these are:
A.~M.~Legendre (1830) \cite[p.~14]{Legendre1830},
E.~Hinkley (1853) \cite[p.~7]{Hinkley1853},
M.~Glaisher (1876) \cite[p.~232]{Lucas1878},
K.~Weierstrass (1876) \cite[p.~391]{Weierstrass1902},
R.~Frick \& F.~Klein (1897) \cite[p.~609]{FK1897},
A.~Cayley (1890) \cite[p.~615]{Cayley1890},
L.~Kronecker (1901) \cite[p.~303]{Kronecker1901},
G.~Chrystal (1904) \cite[p.~38]{Chrystal1904} and
D.~N.~Lehmer (1914) \cite{lehmer1914}.
But on the other hand, we also have quite a number who omitted unity from the primes:
G.~S.~Kl\"ugel (1808) \cite[p.~892]{KGM1808},
P.~Barlow (1811) \cite[p.~54]{Barlow1811},
M.~Ohm  (1834) \cite[p.~140]{Ohm1834},
A.~Reynaud (1835) \cite[pp.~48--49]{Reynaud1835},
L.~Dirichlet (1863) \cite[p.~12]{Dirichlet1863},
E.~Meissel (1870) \cite[p.~34]{ Glaisher1879},
P.~Chebyshev (1889) \cite[pp.~3--4]{Chebyshev1889},
E.~Lucas (1891) \cite[pp.~350--351]{Lucas1891},
J.~P.~Gram (1893) \cite{Gram1893},
E.~Landau (1909) \cite[p.~3]{Landau1909} and
H.~v.~Mangoldt (1912) \cite[p.~176]{Mangoldt1912}.
Together these lists provide a strong argument against the contemporary legend that ``one used to be a prime number''---the truth is far messier.

\section{The Encyclop{\ae}dia Britannica 1890--1910}\label{sect EB}

Another way we might assess the view of the primality of the number one near 1900 is by using the Encyclop{\ae}dia Britannica.  The Encyclop{\ae}dia Britannica's 9th edition, published one volume at a time from 1875 to 1889, is often called the scholar's edition because of its authenticity and originality \cite[p.~62]{Kogan1958}.  The 10th edition (1902--1903) was just a 10--volume supplement to the 25--volume 9th edition, but the encyclopedia was rewritten for the 11th edition (1910--1911).   We will look at the 9th and 11th editions and how they reflected the confusion of the 19th century.

In the 9th edition, Cambridge professor Arthur Cayley wrote an eleven--page article on \textit{number} \cite{Cayley1890} which `defined' prime number in a way which appears (at first) to omit unity:
\begin{quote}
   A number such as 2, 3, 5, 7, 11, \&c., which is not a product of numbers, is said to be a prime number; and a number which is not prime is said to be composite.  A number other than zero is thus either prime or composite.
\end{quote}
After discussing ``complex theories'' (number rings) he returns to the ordinary integers and on the same page writes the following.
\begin{quote}
  Some of these, 1, 2, 3, 5, 7, \&c. are prime, others, $4, =2^{2}, 6, =2.3,$ \&c., are composite; and we have the fundamental theorem that a composite number is expressible, and that in one way only, as a product of prime factors, $N = a^{\alpha}b^{\beta}c^{\gamma} \ldots (a, b, c, \ldots$ primes other than 1; $\alpha, \beta, \gamma, \ldots$ positive integers).
\end{quote}
He is pressed to this definition, like the table makers, by his decision that integers are all either prime or composite, even unity.  Lehmer, in the first paragraph of his introduction to his 1914 impressive table of primes \cite{lehmer1914}, gives exactly this reason for calling unity prime.\footnote{
``A prime number is defined as one that is exactly divisible by no other number than itself and unity.  The number 1 itself is to be considered as a prime according to this definition and has been listed as such in the table.  Some mathematicians [a footnote here cites E.~Landau \cite{Landau1909}], however, prefer to exclude unity from the list of primes, thus obtaining a slight simplification in the statement of certain theorems.  The same reasons would apply to exclude the number 2, which is the only even prime, and which appears as an exception in the statement of many theorems also.  The number 1 is certainly not composite in the same sense as the number 6, and if it is ruled out of the list of primes it is necessary to create a particular class for this number alone.''}  However, in this same edition of the encyclopedia, Rev. George M'Arthur's article on \textit{arithmetic} \cite{MArthur1898}, omits unity without apology:
\begin{quote}
   A prime number is a number which no other, except unity, divides without a remainder; as 2, 3, 5, 7, 11, 13, 17, \&c. [\ldots] The prime factors of a number are the prime numbers of which it is the continued product.  Thus, 2, 3, 7 are the prime factors of 42; 2, 2, 3, 5, of 60.
\end{quote}

In the 11th edition, W.~F.~Sheppard's article on \textit{arithmetic} \cite{Sheppard1910} excludes unity at first, but then curiously states his definition as unusual!
\begin{quote}
   A number (other than 1) which has no factor except itself is called a prime number, or, more briefly, a prime.  Thus 2, 3, 5, 7 and 11 are primes, for each of these occurs twice only in the table.  A number (other than 1) which is not a prime number is called a composite number.  [\ldots]  The number 1 is usually included amongst the primes; but, if this is done, the last paragraph [talking about the fundamental theorem of arithmetic] requires modification.
\end{quote}
G.~B.~Mathews, who wrote the 17--page article on \textit{number} \cite[p.~851]{Mathews1910} (also in the 11th edition), defined prime in an ambiguous way with respect to unity, but then continues on, a few lines later, in a way which clearly excludes unity:
\begin{quote}
   A prime number is one which is not exactly divisible by any number except itself and 1; all others are composite.  [\ldots]  Every number may be uniquely expressed as a product of prime factors.  Hence if $n = p^\alpha q^\beta r^\gamma \ldots$ is the representation of any number $n$ as the product of powers of different primes, the divisors of $n$ are the terms of the product $(1+p+p^2 + \ldots + p^\alpha)(1+q+ \ldots + q^\beta) (1+r+\ldots+r^\gamma)\ldots$ their number is $(\alpha + 1)(\beta+1)(\gamma+1)\ldots,$ and their sum is $\Pi(p^{\alpha+1}{-}1){\div}\Pi(p{-}1).$
\end{quote}

So even as mathematicians were slowly unifying around the modern view, the Encyclop{\ae}dia Britannica continued to illustrate the previous century's confusion.

\section{And the last mathematician was\ldots}

At this point it is natural to ask who was the last significant mathematician to define one as a prime.  (If we leave off any restriction, it is easy to find recent authors still starting the primes with one: expositions such as M.~Kraitchik (1953) \cite[p.~78]{Kraitchik1953}, A.~Beiler (1964) \cite[p.~223]{Beiler1964} and C.~Sagan (1997) \cite[p.~76]{Sagan1997}; as well as the 2011 ``Handy Science Answer Book'' \cite[p.~13]{Handy2011} and Andreasen et al. (2012) \cite[p.~342]{ASO2011}.)

\begin{figure}[h!tb]
  \centering
  \includegraphics[width=.9\textwidth]{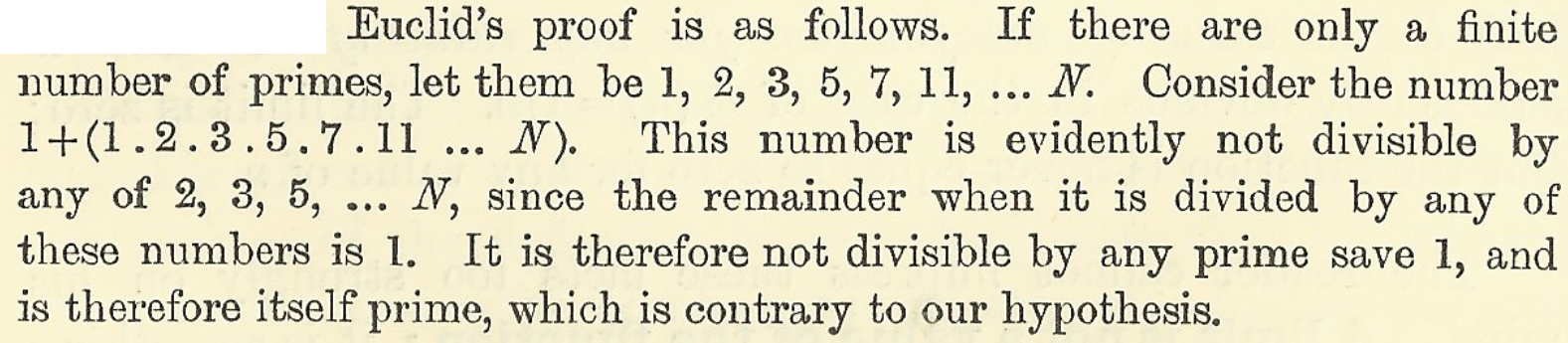}
  \caption{Euclid's proof, in Hardy's \textit{A Course of Pure Mathematics}, sixth edition \cite[p.~120]{Hardy1933}}
  \label{fig_Hardy6}
\end{figure}

\begin{figure}[h!tb]
  \centering
  \includegraphics[width=.9\textwidth]{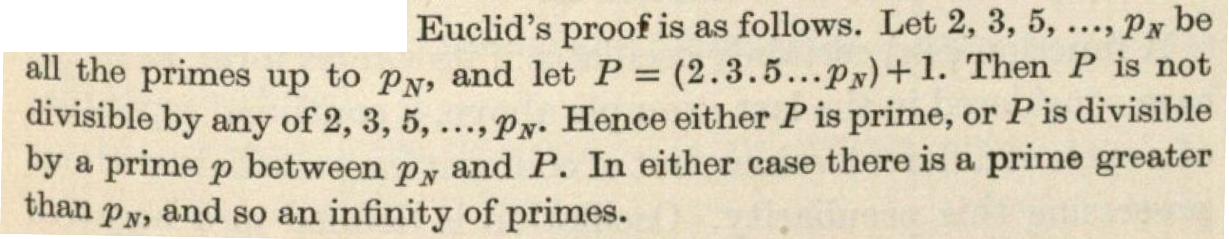}
  \caption{Euclid's proof, in Hardy's \textit{A Course of Pure Mathematics}, seventh edition \cite[p.~125]{Hardy1938}}
  \label{fig_Hardy7}
\end{figure}

There is a legend that the last significant mathematician to define one to be prime was Henri Lebesgue in 1899 (see \cite[p.~33]{Derbyshire2003}).
However, Lebesgue published his first six papers in 1898--1900 \cite{Lebesgue1898,Lebesgue1899a,Lebesgue1899,Lebesgue1899b,Lebesgue1900,Lebesgue1900a}
and none of these contain any references to the prime numbers. On the other hand, well after that time G.\ H.\ Hardy listed one as a prime in
several editions of his text {\it A Course of Pure Mathematics}. For example, in the first six editions (1908, 1914, 1921, 1925, 1928 and 1933), Hardy presented Euclid's proof that there are infinitely many primes with a sequence of primes beginning with 1 (see Figure \ref{fig_Hardy6}).
In the seventh edition (1938 \cite[p.~125]{Hardy1938}), this proof is rewritten so that the primes begin at 2 (see Figure~\ref{fig_Hardy7}).
(He does the same in a 1929 article (\cite[p.~802]{Hardy1929}).)  Hardy (and his editors) missed another reference to the unity as prime, he gives the following example \cite[p.~147]{Hardy1908} of an irrational number: ``the decimal $.111\,010\,100\,010\,10\ldots$, in which the $n$th figure is $1$ if $n$ is prime, and zero otherwise''.  This example remained the same in all ten editions plus the recent ``revised 10th edition'' printed in 2008 \cite[p.~151]{Hardy2008}).

\section{Conclusion}

\begin{figure}[hbt]
  \centering
  \includegraphics[width=.8\textwidth]{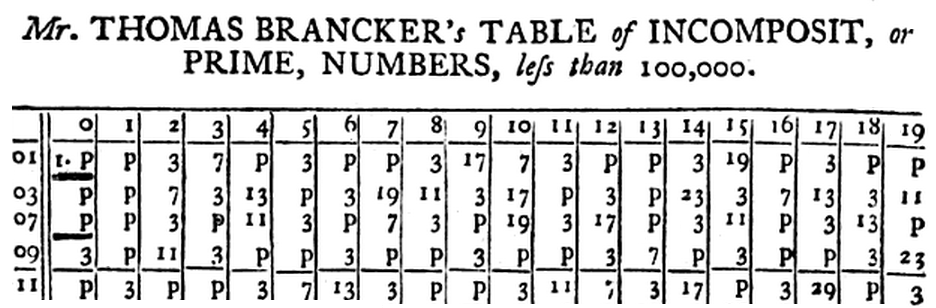}
  \caption{The start of Brancker and Pell's table of factorizations}
  \label{fig_pell}
\end{figure}

There is no question that for the 17th through early 20th century many mathematicians listed the number one as a prime, but it is also clear this definition was never the unified view of mathematicians.  Euclid, Mersenne, Euler, Gauss, Dirichlet, Lucas and Landau all omitted one from the primes, so sequence \seqnum{A008578} of \cite{OEIS2012}, ``prime numbers at the beginning of the 20th century (today 1 is no longer regarded as a prime),'' appears misnamed.  We have seen, though, an obvious choice of a new name for this sequence. Those with the strongest tradition and argument for including unity were the table makers, because what they were really interested in was factorizations, and hence, the numbers that were not composite.  For example, a reprint \cite{Maseres1795} of Brancker and Pell's 1688 table (see Figure \ref{fig_pell}) shows that it is a table of incomposites. Why not call this sequence the incomposites?

When we decided to study the development of the definition of prime (as it applied to finding the smallest prime), we did not expect to find such a deep and rich history.  It seems so obvious now that one is an integer and that two is the first prime, but we use the concepts of integer and prime with an ease of abstraction that took a millennium to develop.

\section{Acknowledgements}
We would like to thank %the referee for his suggestions and
David Broadhurst and Wilfrid Keller for their correspondence on these topics.

%%\begin{thebibliography}{9}%
\RaggedRight
%%%%\bibliographystyle{amsplain}
%%%%\bibliography{one}

\begin{thebibliography}{99}

\bibitem{AF1997}
   A.~G.\ A{\u{g}}arg{\"u}n and C.~R.\ Fletcher, The fundamental theorem of arithmetic dissected, \textit{Math.\ Gazette} \textbf{81} (1997), 53--57.
\bibitem{AO2001}
   A.~G.\ A{\u{g}}arg{\"u}n and E.~M.\ {\"O}zkan, A historical survey of the fundamental theorem of arithmetic, \emph{Historia Math.\ }\textbf{28} (2001), 207--214.
\bibitem{ASO2011}
   J.~B.\ Andreasen, L.~Spalding, and E.~Ortiz, \emph{{FTCE:} {E}lementary {E}ducation {K}-6}, Cliffs Notes, John Wiley \& Sons, 2011.
\bibitem{Barlow1811}
   P.~Barlow, \emph{An {E}lementary {I}nvestigation of the {T}heory of {N}umbers, with its {A}pplication to the {I}ndeterminate and {D}iophantine {A}nalysis,
  the {A}nalytical and {G}eometrical {D}ivision of the {C}ircle, and {S}everal {O}ther {C}urious {A}lgebraical and {A}rithmetical {P}roblems}, J.\ Johnson
  and Co., London, 1811.
\bibitem{Beiler1964}
   A.~Beiler, \emph{Recreations in the {T}heory of {N}umbers}, Dover, 1964.
\bibitem{CRX2012}
   C.~Caldwell, A.~Reddick, and Y.~Xiong, The history of the primality of one---a selection of sources, \url{http://primes.utm.edu/notes/one.pdf}.
\bibitem{Handy2011}
   {Carnegie Library of Pittsburgh}, \emph{The {H}andy {S}cience {A}nswer {B}ook}, Visible Ink Press, 2011.
\bibitem{Cataldi1603}
   P.~A.\ Cataldi, \emph{Trattato de' numeri perfetti}, Heredi di Giovanni Rossi, Bologna, 1603.
\bibitem{Cayley1890}
   A.~Cayley, Number, in \emph{The Encyclop{\ae}dia Britannica}, Vol.~17, The Henry G.~Allen Company, 9th ed., 1890, pp.~614--624.
\bibitem{Chrystal1904}
   G.~Chrystal, \emph{Algebra, an Elementary Text-Book}, Part I, 5th ed., Adam and Charles Black, London, 1904.
\bibitem{CP2005}
   R.~Crandall and C.~Pomerance, \emph{Prime {N}umbers: {A} {C}omputational {P}erspective}, 2nd ed., Springer, 2005.
\bibitem{Derbyshire2003}
   J.~Derbyshire, \emph{Prime Obsession: Bernhard {R}iemann and the Greatest Unsolved Problem in Mathematics}, 2nd ed., Joseph Henry Press, 2003.
\bibitem{Dirichlet1863}
   P.~G.~L.\ Dirichlet, \emph{Vorlesungen {\"u}ber {Z}ahlentheorie}, Friedrich Vieweg und Sohn, Braunschweig, 1863.
\bibitem{Euler1770a}
   L.~Euler, \emph{Elements of {A}lgebra}, 2nd ed., Vol.~1, J.\ Johnson and Company, 1810, translated from the {F}rench; with additions of {L}a {G}range and the notes of the {F}rench Translator.
\bibitem{Felkel1776}
   A.~Felkel, \emph{Tafel aller einfachen {F}actoren der durch 2, 3, 5 nicht theilbaren {Z}ahlen von 1 bis 10 000 000}, Vol.~1, von Ghelen, Wien, 1776.
\bibitem{FK1897}
   R.~Fricke and F.~Klein, \emph{Vorlesungen {\"u}ber die {T}heorie der automorphen {F}unctionen}, Vol.~1, B.~G.\ Teubner, Leipzig, 1897.
\bibitem{Glaisher1879}
   J.~Glaisher, \emph{Factor {T}able for the {F}ourth {M}illion}, Taylor and Francis, London, 1879.
\bibitem{Goldbach1742}
   C.~Goldbach, \emph{Letter to {E}uler}, Collected by P.~H.\ Fuss in his Correspondance math{\'e}matique et physique de quelques c{\'e}l{\`e}bres g{\'e}om{\`e}tres du XVIII{\`e}me si{\`e}cle, 7 June 1742.
\bibitem{Gram1893}
   J.~P.~Gram, \emph{Rapport sur quelques calculs entrepris par {M}.~{B}ertelsen et concernant les nombres premiers}, Acta Math.\ \textbf{17} (1893), 301--314.
\bibitem{Grant1974}
   E.~Grant, \emph{A {S}ource {B}ook in {M}edieval {S}cience}, Harvard Univ.\  Press, 1974.
\bibitem{Hardy1908}
   G.~H.\ Hardy, \emph{A {C}ourse of {P}ure {M}athematics}, 1st ed., Cambridge Univ.\ Press, 1908.
\bibitem{Hardy1929}
   \bysame, An introduction to number theory, \emph{Bull.\ Amer.\ Math.\ Soc.} \textbf{35} (1929), 778--818.
\bibitem{Hardy1933}
   \bysame, \emph{A {C}ourse of {P}ure {M}athematics}, 6th ed., Cambridge Univ.\ Press, 1933.
\bibitem{Hardy1938}
   \bysame, \emph{A {C}ourse of {P}ure {M}athematics}, 7th ed., Cambridge Univ.\ Press, 1938.
\bibitem{Hardy2008}
   \bysame, \emph{A {C}ourse of {P}ure {M}athematics}, revised 10th ed.\ with a foreword by T.~W.\ K\"orner, Cambridge Univ.\ Press, 2008.
\bibitem{Heath1981}
   T.~L.\ Heath, \emph{A {H}istory of {G}reek {M}athematics: From {T}hales to {E}uclid}, Vol.~1, Dover, 1981.
\bibitem{Hinkley1853}
   E.~Hinkley, \emph{Tables of the {P}rime {N}umbers, and {P}rime {F}actors of {C}omposite {N}umbers, from 1 to 100,000: {W}ith the {M}ethods of {T}heir
  {C}onstruction, and {E}xamples of {T}heir {U}se},  Baltimore, 1853.
\bibitem{Jones1978}
   C.~J.~Jones, \emph{The concept of \underline{one} as a number}, Ph.D.\ thesis, Univ.\ Toronto, 1978.
\bibitem{KK2011}
   K.~U.~Katz and M.~G.~Katz, Stevin numbers and reality, \emph{Found.\ Sci.} \textbf{17}, 2012.
\bibitem{AlKindi1974}
   Kind{\=\i} and A.~L.\ Ivry, Al-{K}ind{\=\i}'s {M}etaphysics: {A} {T}ranslation of {Y}a`qub ibn {I}sa\d{h}{\=a}q al-{K}ind{\=\i}'s {T}reatise
  ``{O}n {F}irst {P}hilosophy'' (f{\=\i} al-{F}alsafah al-{\=u}l{\=a}), in \emph{Studies in Islamic Philosophy and Science}, State Univ. N.Y.
  Press, 1974.
\bibitem{KGM1808}
   G.~S.\ Kl{\"u}gel, \emph{Mathematisches {W}{\"o}rterbuch}, Part I: \emph{Die reine Mathematik}, Vol.~3, E.~B.\ Schwickert, Leipzig, 1808.
\bibitem{Kogan1958}
   H.~Kogan, \emph{The Great {EB}: The Story of the {E}ncyclopaedia {B}ritannica}, Univ. of Chicago Press, 1958.
\bibitem{Kraitchik1953}
   M.~Kraitchik, \emph{Mathematical {R}ecreations}, Dover, New York, 1953.
\bibitem{Kronecker1901}
   L.~Kronecker, \emph{Vorlesungen {\"u}ber {M}athematik}, Part II: \emph{Allgemeine Arithmetik}, Vol.~1, B.~G.\ Teubner, Leipzig, 1901.
\bibitem{Lambert1770}
   J.~H.\ Lambert, \emph{Zus{\"a}tze zu den logarithmischen und trigonometrischen {T}abellen}, Haude und Spener, Berlin, 1770.
\bibitem{Landau1909}
   E.~Landau, \emph{Handbuch der {L}ehre von der {V}erteilung der {P}rimzahlen}, Vol.~1, B.~G.\ Teubner, Leipzig and Berlin, 1909.
\bibitem{Lebesgue1898}
   H.~Lebesgue, Sur l'approximation des fonctions, \emph{Bull. Sci. Math.} \textbf{22} (1898), 278--287.
\bibitem{Lebesgue1899a}
   \bysame, Sur la d\'efinition de l'aire d'une surface, \emph{C.\ R.\ Math.\ Acad.\ Sci.\ Paris}  \textbf{129} (1899), 870--873.
\bibitem{Lebesgue1899}
   \bysame, Sur les fonctions de plusieurs variables, \emph{C.\ R.\ Math.\ Acad.\ Sci.\ Paris}  \textbf{128} (1899), 811--813.
\bibitem{Lebesgue1899b}
   \bysame, Sur quelques surfaces non r\'egl\'ees applicables sur le plan, \emph{C.\ R.\ Math.\ Acad.\ Sci.\ Paris} \textbf{128} (1899), 1502--1505.
\bibitem{Lebesgue1900}
   \bysame, Sur la d\'efinition de certaines int\'egrales de surface, \emph{C.\ R.\ Math.\ Acad.\ Sci.\ Paris} \textbf{131} (1900), 867--870.
\bibitem{Lebesgue1900a}
   \bysame, Sur le minimum de certaines int\'egrales, \emph{C.\ R.\ Math.\ Acad.\ Sci.\ Paris} \textbf{131} (1900), 935--937.
%\bibitem{Lebesgue1972}   % removed statement the above six are more easily found in these two
%   \bysame, \emph{{\OE}uvres scientifiques}, Vol.~3, L'Enseignement Math\'ematique, 1972.
%\bibitem{Lebesgue1973}
%   \bysame, \emph{{\OE}uvres scientifiques}, Vol.~4, L'Enseignement Math\'ematique, 1973.
\bibitem{Lebesgue1856}
   V.~A.\ Lebesgue, Remarques diverses sur les nombres premiers, \emph{Nouv. Ann. Math.} \textbf{15} (1856), 130--143.
\bibitem{Lebesgue1859}
   \bysame, \emph{Exercices d'analyse num{\'e}rique}, Leiber et Faraguet, 1859.
\bibitem{Lebesgue1864}
   V.~A.\ Lebesgue and J.~Ho{\"u}el, \emph{Tables diverses pour la d{\'e}composition des nombres en leurs facteurs premiers}, Gauthier-Villars, 1864.
\bibitem{Legendre1830}
   A.~M.\ Legendre, \emph{Th{\'e}orie des nombres}, 3rd ed., Vol.~1, Firmin {D}idot Fr{\`e}res, Paris, 1830.
\bibitem{lehmer1914}
   D.~N.\ Lehmer, \emph{List of {P}rime {N}umbers from 1 to 10,006,721}, Carnegie Institution of Washington, 1914.
\bibitem{Lucas1878}
   E.~Lucas, Th\'eorie des fonctions num\'eriques simplement p\'eriodiques, \emph{Amer. J. Math.} \textbf{1} (1878), 197--240.
\bibitem{Lucas1891}
   \bysame, \emph{Th\'eorie des nombres}, Gauthier-Villars, Paris, 1891 (reprinted: 1991 Jacques Gabay, Sceaux).
\bibitem{Mangoldt1912}
   H.~v.~Mangoldt, \emph{Einf{\"u}hrung in die h{\"o}here {M}athematik}, Vol.~2: \emph{{D}ifferentialrechnung}, S.~Hirzel, Leipzig, 1912.
\bibitem{MArthur1898}
   G.~R.~M'Arthur, Arithmetic, in \emph{The Encyclop{\ae}dia Britannica}, Vol.~2, A.~and C.~Black, 9th ed., 1898, pp.\ 524--536.
\bibitem{Maseres1795}
   F.~Maseres, J.~Bernoulli, and J.~Wallis, \emph{The {D}octrine of {P}ermutations and {C}ombinations, {B}eing an {E}ssential and {F}undamental {P}art of the {D}octrine of {C}hances}, Francis Maseres, London, 1795.
\bibitem{Masi1983}
   M.~Masi, \emph{Boethian {N}umber {T}heory---{A} {T}ranslation of the {D}e {I}nstitutione {A}ritmetica}, Studies in Classical Antiquity, Vol.~6, Rodopi,
   New York, 1983.
\bibitem{Mathews1910}
   G.~B.\ Mathews, Number, in \emph{The Encyclop{\ae}dia Britannica}, Vol.~19, Cambridge Univ. Press, 11th ed., 1910, pp.~847--866.
\bibitem{Menninger1992}
   K.~Menninger, \emph{Number {W}ords and {N}umber {S}ymbols: {A} {C}ultural {H}istory of {N}umbers}, Dover, 1992.
\bibitem{Mersenne1625}
   M.~Mersenne, \emph{La {V}{\'e}rit{\'e} des sciences: contre les sceptiques ou pyrrhoniens}, Toussaint du Bray, Paris, 1625.
\bibitem{Moxon1679}
   J.~Moxon, \emph{Mathematicks {M}ade {E}asy, or, a {M}athematical {D}ictionary {E}xplaining the {T}erms of {A}rt and {D}ifficult {P}hrases {U}sed in {A}rithmetick, {G}eometry, {A}stronomy, {A}strology, and {O}ther {M}athematical {S}ciences\ldots}, Joseph Moxon, 1679.
\bibitem{Ohm1834}
   M.~Ohm, \emph{Die reine {E}lementar-{M}athematik}, Vol.~1: \textit{Die Arithmetik bis zu den h\"ohern Gleichungen}, 2nd ed., Jonas Verlags-Buchhandlung, Berlin, 1834.
\bibitem{Prestet1689}
   J.~Prestet, \emph{Nouveaux elemens des mathematiques}, Andr\'e Pralard, Paris, 1689.
\bibitem{RX2012}
   A.~Reddick and Y.~Xiong, The search for one as a prime number: from ancient {G}reece to modern times, preprint.
\bibitem{Reynaud1835}
   A.~A.~L.\ Reynaud, \emph{Trait{\'e} d'arithm{\'e}tique}, Bachelier, Paris, 1835.
\bibitem{Roegel2011}
   D.~Roegel, \emph{A Reconstruction of the Tables of the {S}huli {J}ingyun (1713-1723)}, Research report, LORIA, 2011, \url{http://hal.inria.fr/hal-00654450/en}.
\bibitem{Sagan1997}
   C.~Sagan, \emph{Contact}, Pocket Books, 1997.
\bibitem{Leon1657}
  L.~de Saint-Jean, \emph{Studium sapienti{\ae} universalis}, Iacobus Quesnel, Paris, 1657.
\bibitem{Schooten1657}
   F.~v.~Schooten, \emph{Exercitationes mathematicae}, J.~Elsevirius, Leiden, 1657.
\bibitem{Sheppard1910}
   W.~F.\ Sheppard, Arithmetic, in \emph{The Encyclop{\ae}dia Britannica}, Vol.~2, Cambridge Univ. Press, 11th ed., 1910, pp.~523--542.
\bibitem{OEIS2012}
   N.~J.~A.\ Sloane, \emph{The online encyclopedia of integer sequences}, electronically published at \url{http://oeis.org}.
\bibitem{Smith1958}
   D.~E.\ Smith, \emph{History of {M}athematics}, Vol.~II, Dover, 1958.
\bibitem{Stahl1992}
   W.~H.\ Stahl and E.~L.\ Burge, \emph{Martianus {C}apella and the {S}even {L}iberal {A}rts}, Vol.\ II, Columbia Univ. Press, 1992.
\bibitem{Taran1981}
   L.~Tar{\'a}n, \emph{Speusippus of {A}thens: {A} {C}ritical {S}tudy with a {C}ollection of the {R}elated {T}exts and {C}ommentary}, Philosophia Antiqua,
  Vol.~39, E.~J.~Brill, Leiden, 1981.
\bibitem{Chebyshev1889}
   P.~L.\ Tschebyscheff, \emph{Theorie der {C}ongruenzen}, edited and translated by H.~H.\ Schapira, Mayer \& M{\"u}ller, Berlin, 1889.
\bibitem{Waring1782}
   E.~Waring, \emph{Meditationes algebraic{\ae}}, J.\ Archdeacon, Cambridge, 1782.
\bibitem{Weierstrass1902}
   K.~Weierstrass, \emph{Vorlesungen \"uber die Theorie der Abelschen Transcendenten}, Math. {W}erke, Vol.~4, Berlin, 1902.
\bibitem{OldWiki}
   Wikipedia contibutors, Prime number, \emph{{W}ikipedia{,} The Free Encyclopedia}, 2011.
\end{thebibliography}
\providecommand{\bysame}{\leavevmode\hbox to3em{\hrulefill}\thinspace}
\providecommand{\MR}{\relax\ifhmode\unskip\space\fi MR }
\providecommand{\MRhref}[2]{%
  \href{http://www.ams.org/mathscinet-getitem?mr=#1}{#2}
}
\providecommand{\href}[2]{#2}

\bigskip
\hrule
\bigskip

\noindent 2010 {\it Mathematics Subject Classification}:
Primary 11A41; Secondary 11A51, 01A55.

\noindent \emph{Keywords: } prime number, unity, history of mathematics.

\end{document}